\documentclass[reprint, frontmatterverbose, nofootinbib, amsmath, amssymb, aps]{revtex4-1}

\usepackage{graphicx}
\usepackage{dcolumn}
\usepackage{bm}
\usepackage[mathlines]{lineno}
\usepackage{lipsum}
\usepackage{todonotes}
\usepackage{bbding}

\def\cc{\mathbf{c}}
\def\me{\textit{median }}
\def\de{\delta}

\newcommand{\R}{{\mathbb R}}

\begin{document}
\title{Cost-Optimal Operation of Energy Storage Units:\\Impact of Uncertainties and Robust Estimator}

\author{Lars Siemer}
\email{Corresponding author: lars.siemer@uni-bremen.de}
\affiliation{%
 Departement of Mathematics, Univerity of Bremen, 28334 Bremen, Germany.
}
\author{Wided Medjroubi}%
\affiliation{DLR Institute of Networked Energy Systems, Department Energy Systems Analysis, Carl-von-Ossietzky-Str.~15, 26129 Oldenburg, Germany}

\begin{abstract}
The rapid expansion of wind and solar energy leads to an increasing volatility in the electricity generation. Previous studies have shown that storage devices provide an opportunity to balance fluctuations in the power grid. An economical operation of these devices is linked to solutions of probabilistic optimization problems, due to the fact that future generation is not deterministic in general. For this reason, reliable forecast methods as well as appropriate robust optimization algorithms take an increasingly important role in future power operation systems. Taking an uncertain availability of electricity into account, we present a method to calculate cost-optimal charging strategies for energy storage units. The proposed method solves subproblems which result from a sliding window approach applied on a linear program by utilizing statistical measures. The prerequisite of this method is a recently developed fast algorithm for storage-related optimization problems. To present the potential of the proposed method, a Power-To-Heat storage system is investigated as an example using today's available forecast data and a robust statistical measure. Second, the novel approach proposed here is compared with a common robust optimization method for stochastic scenario problems. In comparison, the proposed method provides lower acquisition costs, especially for today's available forecasts, and is more robust under perturbations in terms of deteriorating predictions, both based on empirical analyses. Furthermore, the introduced approach is applicable to general cost-optimal operation problems and also real-time optimization concerning uncertain acquisition costs. 
\end{abstract}

\maketitle

\section{\label{sec1}Introduction}
Electricity generation becomes more and more unpredictable in future, due to increasing shares of Renewable Energy Sources (RES) mainly from wind and solar resources, in the electricity generation portfolio. Additionally, RES generation, due to its fluctuating and non-predictable nature, is not adaptable to temporally changes in the electricity demand. The integration process of high shares of RES is still ongoing which is making options pertaining to their economical operation management important~\cite{lund2007renewable, georgilakis2008technical, boyle1997renewable, monteiro2009wind}. 

One option available to utilize overproduction and to prevent curtailment are energy storage units. They can contribute substantially to the expansion of RES within the electricity grid and increasingly required to maintain grid stability in prospective power supply systems~\cite{braff2016value, weitemeyer2015integration, weitemeyer2017optimal}. Because of the volatile nature of RES, robust forecast as well as suitable optimization methods are essential in the future with respect to an economical operation of these storage devices~\cite{jones2014renewable}.

The focus will be on suitable optimization methods, hence on methods which calculate robust and cost-optimal solutions within a reasonable time frame. Well-known approaches are for example stochastic optimization methods or the Min-Max method. The latter one minimizes the worst-case scenario and therefore provides the most robust solution with respect to perturbations in the initial data~\cite{schneider2007stochastic, ben2009robust}.

In this Letter, a novel approach is described which is based on a generalization of Model Predictive Control (MPC) approaches~\cite{grune2011nonlinear}. In detail, a robust statistical measure is applied to subsolutions of a truncated optimization problem, resulting from a sliding window approach on all trajectories of simulated price realizations. This approach leads to an increase in the number of calculations, especially for high temporal resolutions. For reasonable computation times, a recent developed fast optimization method for storage related problems is utilized here. This method solves large problems in a short computational time, due to a quadratic upper bound in runtime~\cite{siemer2016jes}.

As a particular example for the application of energy storage units to integrate large shares of RES efficiently, especially in times of overproduction, Power-to-Heat (PtH) storage units installed in private households are considered here. Assuming in addition that the electricity price reflects the fluctuating availability and unpredictability of energy, these systems need to be operated in a  cost-optimal way. Annual acquisition costs for the considered system are calculated in a short time by the novel approach utilizing the median as a robust statistical measure. To characterize uncertain predictions of electricity prices, calculations are conducted which are based on a large number of portfolios for electricity prices at each point in time, simulated using  commercial  forecast  data. Additionally, the results are further compared with the Min-Max method in terms of lower overall acquisition costs as well as robustness with respect to uncertain price forecasts. The novel approach leads to lower costs while more robust solutions, both with respect to today common forecast data.

Beyond this particular field of application, the proposed method is applicable to general storage related problems and describes a novel approach to solve stochastic predictive control problems efficiently, based on uncertain forecast trajectories.

\section{\label{sec2}Methodology}

In this Letter, we study the dynamics of solutions of linear optimization problems (\ref{LP}), including uncertainties in the objective function. In particular, the focus lies on the analysis of cost-optimal charging strategies for energy storage units using a new method. For this reason, we focus on the following type of minimization problems

\begin{align}\label{LP}
\min \quad & \sum_{i=1}^{n}{c_i \cdot x_i}   \nonumber \\
\mbox{such that} \quad & 0 \le x_i \le u_i\,, \ \mbox{$\forall \, i=1, \dots, n$} \tag{LP}\\
\mbox{and}\quad  & a_i \le \sum_{j=1}^{i}{q^{i-j} \cdot x_j} \le b_i\,, \ \mbox{$\forall \, i=1, \dots, n$}  \nonumber
\end{align}

where, $n$ is the length of the equally discretized period of time considered, $c \in \R^n$ are the acquisition costs, and  $a\,,b\,,u \in \R^n$, as well as $0<q \le 1$ define technical restrictions for the energy storage unit. The vector $x \in \R^n$ defines the amount of energy used to charge the storage unit for each time-interval (for details, see \cite{siemer2016jes}).

The basic idea of the proposed method is to use a statistical measure on subsolutions of the (\ref{LP}). These subsolutions result from a sliding window approach when considering different portfolios of possible predicted realizations of the vector $c$.

In detail, let $m \in \mathbb{N}$ be the number of simulated forecast scenarios at an initial point in time, e.g.~$m$ different future spot market prices, and $r \in \mathbb{N}\,, r \le n$, a (not necessarily fixed) forecast period, e.g.~representing one week, where $n$ is the overall observed period of time. Thus, a portfolio of predicted prices for a given point $k$ in time is defined by $\cc ^k:=\left( c^{k,1}, \dots , c^{k,m} \right) \in \R^r$. Notice that, the considered time period needs to be extended by the forecast period $r$ as well as the constraints. Additionally, it is important to mention that the sliding window size should not be larger as the provided forecast window of the underlying predicted data set in order to avoid biased results. The proposed method is now defined as follows: 

First, set $k=1$ and solve the corresponding subproblem of (\ref{LP}) for $\cc^1$ and $i=1,\dots,r$ of dimension $r$ for each scenario $m$. This leads to $m$ solutions $x^{k,j} \in \mathbb{R}^r$ for $j=1,\dots,m$. Second, apply a statistical measure $\Phi$ on $x^{k,j}_1$ and set $x^\ast_k:=\Phi(x^{k,j}_1)$, where $x^\ast \in \mathbb{R}^n$ represents the final approximated solution. Third, for the sliding window step, set $k \mapsto k+1$ and set $i=k+1,\dots,r+k+1$ for the next subproblem of (\ref{LP}) including an acquisition vector $\cc^{k+1}$. Repeat the previous steps until $k=n+1$. This approach leads to a solution
\[x^\ast:= \left( \Phi(x_1^{1,j}), \dots, \Phi(x_1^{n,j}) \right), \quad j=1,\dots ,m\] 
for the observed period $n$. Due to the linearity of the problem and the convexity of the constraints, the approximated solution $x^\ast$ is a feasible solution of (\ref{LP}). It is important to mention that, the constraints in (\ref{LP}) change for each $k$. Furthermore, this approach implies that, the number of problems to be solved increases tremendously.

The authors recently proposed a new optimization algorithm, which solves problems such as (\ref{LP}) even for large $n$ in a fraction of time required by common solvers~\cite{siemer2016jes}. Therefore, the method proposed above in combination with the optimization algorithm can be realized in a short computational time. This is true even for a large portfolio of simulated scenarios and a long forecast period, which with conventional optimization algorithms would not be feasible within a reasonable period of time.

\section{\label{sec3}Results}

In this section, we first study the impact of predicted electricity spot market prices on the cost-optimal operation through the application of the previously proposed method. The basis for the calculations is a mathematical model of a Power-To-Heat storage unit (PtH), which will serve as an example for a suitable energy storage system. Second, the obtained results will be compared with results of a well-known scenario approach for stochastic optimization problems.

Following~\cite{siemer2016jes}, we study the dynamics of the (\ref{LP}), where $x^\ast \in \R^n$ now describes the resulting charge strategy of the storage in [kW] and $c \in \R^n$ represents the hourly electricity acquisition costs in [EUR/kWh] for a period of time discretized into $n$ intervals of equal size. In addition, $u\in \R^n$ defines the maximum charge power [kW] of the PtH conversion unit and is here set constant. The vectors $a,b \in \R^n$ represent the technical limitations of the energy storage unit [kWh] and $0<q\le 1$ is the corresponding hourly energy reduction factor of the related storage unit.

\begin{figure}[h]
    \centering
    \includegraphics[width=0.48\textwidth]{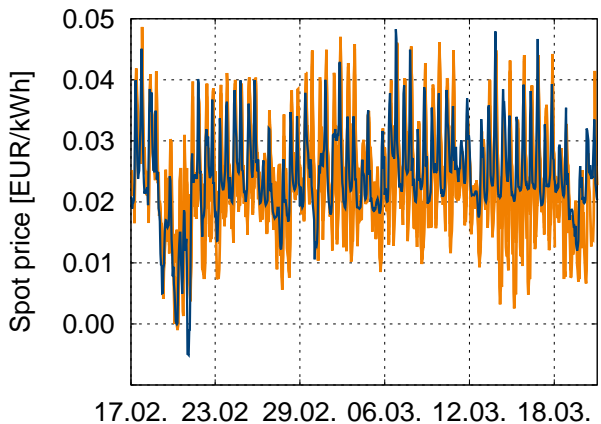}
    \vspace*{-7mm}
    \caption{Hourly price development of the day-ahead electricity price (PHELIX) at the spot market EPEX SPOT (blue) and corresponding daily forecasts for 5-days(orange), (covering the period 17.02.2016 – 23.03.2016).}
    \label{fig:forecast}
\end{figure}

As a data basis for the electricity acquisition costs, we consider the Physical Electricity Index (PHELIX) which is the averaged value of the day-ahead shares traded on the EPEX Spot Market and used as a proxy to represents the volatility of the RES feed-in. In order to study the influences of uncertainties, we first analyse the properties of a commercial PHELIX spot price prediction. For this reason, a data set from 18.02.2016 till 16.03.2016, which represents the daily forecasts for five days in the future, is used (provided by the company Energy Brainpool\footnote{Energy Brainpool GmbH $\&$ Co.~KG., \emph{http://www.energybrainpool.com}}). This leads to a prediction horizon of 3480 hours for the provided data set and this predictions are used to characterize uncertainties on the electricity acquisition costs in the following (cf.~figure \ref{fig:forecast}).

The prediction-dynamics of the given forecast data set are modelled by a stochastic differential equation, characterized by a linear repulsive drift and a constant diffusion term, as in~\cite{kleinhans2012estimation, doob1942brownian}. This leads to a (discrete) mean-reverting Ornstein-Uhlenbeck (OU) stochastic process, including a mean-reverting level and also an initial value of zero~\cite{gillespie1996exact}. In addition, the time increments are equally set to one due to the hourly scaling of the provided data set. The drift and diffusion terms can be analytically derived, where the formula is given by equation (7) in~\cite{kleinhans2012estimation}. We applied the formula to the differences of the predicted spot market prices and the actual ones, for each hour of the 3480 hours given by the provided forecast set mentioned above.

The further averaging of these values leads to a drift of $0.3331$ and a effective diffusion of about $0.004$ for the considered time period. These values are used here to include uncertainties in the objective functions of the (\ref{LP}) subproblems, by using known spot market prices, simulating an OU process in terms of the values above, and adding it to the original price signal. For this reason, the initial value and the mean reversion level of the stochastic process is zero at the beginning. 

To demonstrate the potential with respect to cost-optimal solutions and to observe their dynamics, the introduced method is applied to the following energy storage system. The example system consists of a PtH storage system installed in a typical private single-family household with two persons, living space of 100~m$^2$, specific heat demand of 100~kWh/(m$^2 \cdot$a), and annual heat demand of 10741.15~kWh/a over a period of one year (hourly discretization) from 01.07.2012 to 30.06.2013~\cite{siemer2016jes}. The considered system is composed of a charge unit with maximum charge power of 15~kW, an attached thermal storage unit with a maximum capacity of 3-day mean-average daily heat demand of the considered household ($\approx 88.28$~kWh), and a corresponding hourly energy reduction factor of $q=0.9981$. This system is chosen, because it features a large cost-saving potential while being technically simple to implement (see~\cite{siemer2016jes} and the references therein for details). It should be noted that, the efficiency of the considered conversion unit is set here to 100$\%$ for simplicity reasons. This can be achieved by utilizing for example a heat pump system with a high performance factor~\cite{siemer2016jes}. Additionally, the number of predicted portfolios is fixed to $m=1000$ in order to study a representative sample. The forecast-window for each time-step is further set (constant) to five days in this work, which is equal to the forecast window of the spot price prediction data sets discussed above. As a statistical measure, the \emph{median} is used to represent a robust measure with a breakdown point of 50~$\%$~\cite{maronna2006robust}.

The discrete OU process is generated for each prediction $m$ according to~\cite{gillespie1996exact} and is added as a multiple to the realized PHELIX price signal of the considered time period. In detail, the simulated stochastic processes are linearly added as multiples of the value $\de \in \left[0.1\,,2.5\right]$ with a discrete step-size of $0.1$ in order to facilitate comparison, where $\de=1$ can be considered as today's forecast accuracy resulting from the provided data. $\de=0$ reflects a known price signal in advance, which leads to acquisition costs of 215.62~EUR and defines the limit for $\de \rightarrow 0$. The other extreme is described by $\de \rightarrow \infty$ which leads to acquisition costs of around 443.13~EUR and corresponds to a system without any energy storage unit.

For comparison, the Min-Max method, which is a well-known method to minimize the possible loss for a worst case scenario, is also used in this letter. This method minimizes the maximum value by considering the simulated prices for each $\de$ value as described above~\cite{aissi2009min}. In detail, the Min-Max approach is applied to the model to obtain a robust value of the electricity acquisition costs for further comparison.

\begin{figure}
    \centering
    \includegraphics[width=0.48\textwidth]{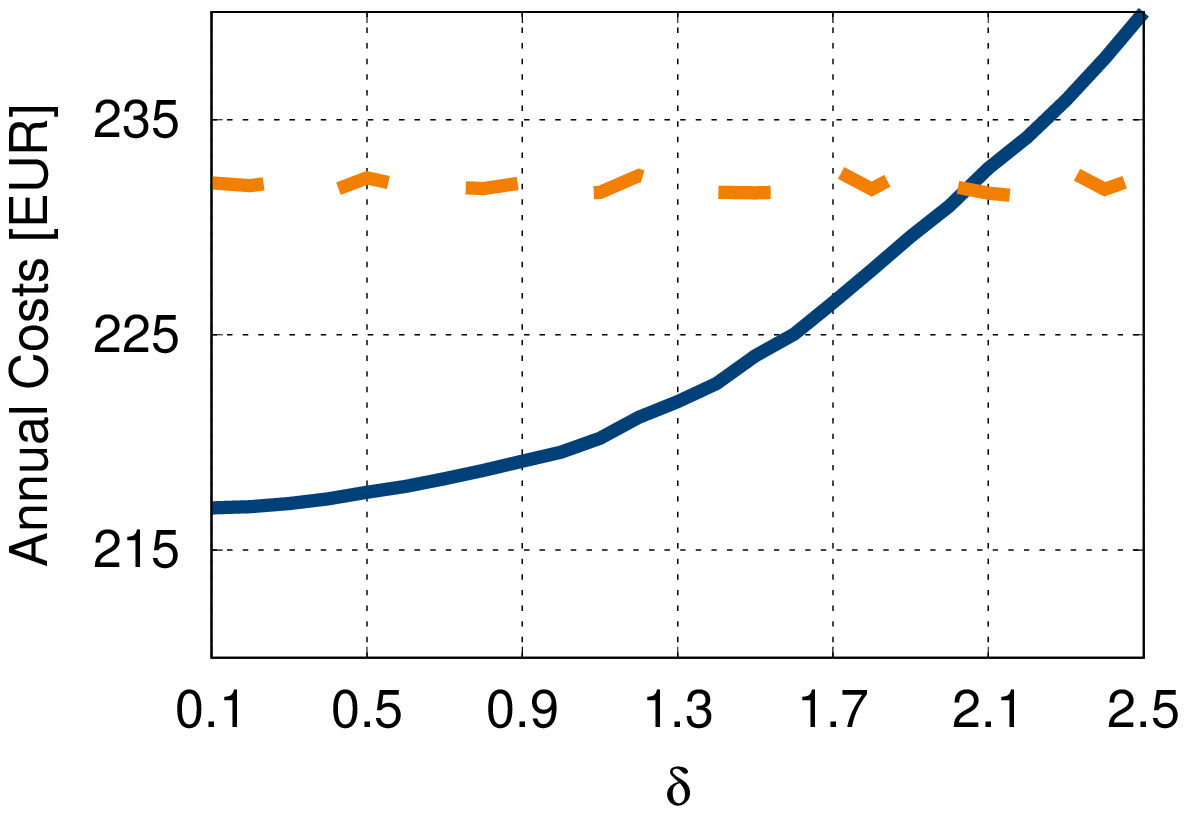}
    \vspace*{-7mm}
    \caption{Approximated annual costs as solutions of the \me (solid blue) and the Min-Max approach (dashed orange), both depending on the parameter $\de$. \\$[$Assumptions: PtH storage unit including max.~charge power of 15~kW, storage capacity of three days (88.28~kWh), an energy storage loss factor of 0.9981, a prediction horizon of five days, and a simulation number of 1000 prices.$]$}
    \vspace*{-4mm}
    \label{fig:result}
\end{figure}

The approximated charging strategies as results of the \emph{median} and of the Min-Max approach are further multiplied by the occurring price signal of the considered period of time (see figure \ref{fig:result}). It should be noted that, the application of a sliding-window approach may lead to a storage level above zero at the end of the considered time period. However, this surplus of energy was not taken into account, due to the low amounts of about $10^{-5}$~EUR for the Min-Max method and about 0.8084~EUR for the novel approach. Both values are based on the mean of the PHELIX price.

\vspace*{-5mm}
\begin{figure}
    \centering
    \includegraphics[width=0.48\textwidth]{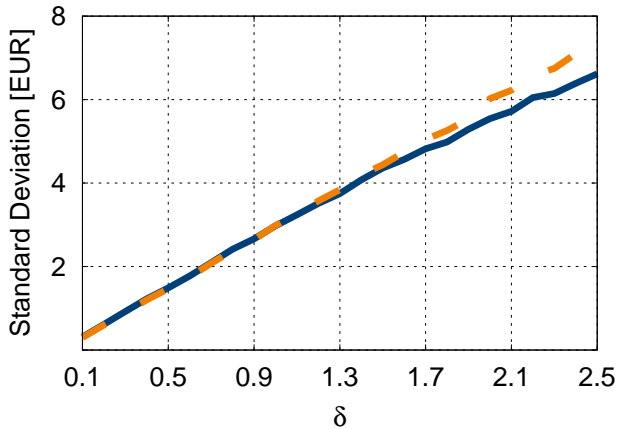}
    \vspace*{-7mm}
    \caption{Standard deviation of the annual costs for the solutions of \me (solid blue) and the Min-Max method (dotted orange) depending on the value of $\de$. Both are based on a sample of five thousand different price realizations for each $\de$. Further assumptions are analogous to figure \ref{fig:result}.}
    \vspace*{-4mm}
    \label{fig:deviation}
\end{figure}
\newpage
The curves in figure \ref{fig:result} describe the annual electricity acquisition costs and are resulting from the application of the \me and the Min-Max method. For the method proposed here, the curve tends to the value 215.62~EUR for $\de \rightarrow 0$. For increasing value of $\de$, the price is about 443.13~EUR, corresponding to bigger uncertainties in the predicted price signal. In contrast, the results of the Min-Max approach are almost constant for all values of $\de \in [0.1,2.5]$ at a price of about 231.96~EUR. This can be explained by the worst-case characteristic of the Min-Max method as well as the selected window-size of five days. Furthermore, the results of Min-Max method tend to the upper limit for an increasing $\de$, but with a less steeper slope compared with the \emph{median} method. Additionally, both results intersect for $\de$ value between $2.0$ and $2.1$. Hence, the proposed method leads to lower annual costs compared with the Min-Max method in this interval. In conclusion, for the lowest level of uncertainties at $\de=0.1$, the \me leads to about $7\,\%$ less overall annual costs and even considering today's forecast accuracy $\de=1.0$, the application of the new method still results in $5.5\,\%$ annual cost savings.

To compare in more details the \me and the minimum worst-case case in the context of sensitivity, the standard deviation for different price scenarios is discussed in the following. For this purpose, a portfolio of five thousand different price scenarios was simulated as described above and the standard deviation related to the results of the \me and the Min-Max approach, both depending on $\de$, were calculated. The results are presented in figure \ref{fig:deviation} and show a linear dependency on $\de$ for the Min-Max method.

The standard deviation resulting from the application of \me is concave. Both curves in figure \ref{fig:deviation} tend to zeros in the limit $\de \rightarrow 0$ and are strictly monotonically increasing for increasing values of $\de$. Furthermore, the standard deviation of \me is lower (or equal at value of $\de$ close to zero) compared with the Min-Max method. The \me leads to a lower deviation for value of $\de$ close to today's forecast accuracy and to about $13.5\,\%$ less for an addition of $\de=2.5$, due to the negative curvature. This results in a lower overall sensitivity of the proposed method compared with the worst case approach, especially for deteriorating predictions.

\section{\label{sec4}Conclusions}
In this Letter, we presented a novel approach to solve a class of linear optimization problems, which introduces uncertainties in the objective function, by using a recently developed optimization algorithm for general storage related problems. The proposed method uses statistical measures on solutions of the subproblems resulting from a sliding-window approach on the primary optimization problem. With the focus on robustness, we used the median as a statistical measure.

In order to increase and economically and efficiently use electricity from renewable energy sources, we applied the developed approach to a mathematical model which describes conversion systems combined with thermal storage units as an example. The aim was to analyse the impact of uncertainties in the acquisition costs on the optimal system operation strategies. Therefore, we first assumed that the PHELIX day-ahead spot market price represents the availability of electricity as well as the distinctive volatility of the RES feed-in. The electricity price are not deterministic and need to be predicted. We assumed that, the forecast uncertainty is growing in time but bounded in time within a confidence interval. For this reason, first a mean-reverting Ornstein-Uhlenbeck process with a linear drift and a constant diffusion term was used to reflect the uncertainties and a data set of real PHELIX spot prices forecasts was analysed to adjust the parameter space for the utilized mean-reverting stochastic process. Second, a PtH conversion system in combination with a thermal storage unit installed in a private household was chosen as an example of a technically simple and feasible opportunity for the application of storage devices in the future. The considered time period was set to one year. The statistically robust measure of the median was used within the evaluation step of this novel approach. 

The results of the proposed method indicate a large cost-saving potential for decreasing forecast uncertainties. Furthermore, a strong sublinear behaviour of the sensitivity for different price realizations in terms of deteriorating predictions can be observed.

In a second step, the results of the novel approach were also compared with the Min-Max method, a common method for stochastic optimization problems. This method is based on the idea of calculating a robust and (cost-) minimal worst-case solution. The Min-Max method leads to almost constant results for the annual acquisition costs as well as linear dependency of uncertainties with respect to sensitivity due to the linear addition of the simulated mean-reverting processes and to the considered window size. In comparison, the new method presented in this letter achieves lower annual acquisition costs up to an addition of about two times higher uncertainties within the price signal as the level of today common PHELIX price predictions. This leads to lower overall electricity acquisition costs for definite price predictions by the application of the novel approach as well as the utilization of a robust statistical measure. On the other hand, the sensitivity performance of the proposed method is better than the one for the Min-Max approach in general which indicates a higher robustness of the developed method compared with the worst-case algorithm, especially for an increase of uncertainty in the price prediction.

In conclusion, the novel optimization method developed in this letter leads to lower acquisition costs including today common price forecasts and the results are more robust with respect to perturbations in the predicted prices. This method is applicable to general linear storage-related optimization problems dealing with uncertainties in the objective function. Thus, improvements in price predictions would lead to an economical and also robust operation of storage devices with respect to perturbed acquisition cost forecasts by the application of this novel approach. Additionally, further development of an appropriate problem-related statistical measure used by the proposed method would also increase the performance.
\section{\label{sec5}Acknowledgements}
We gratefully thank David Kleinhans for assistance and inspiring discussions which contributed to the methodology described above as well as for comments which greatly improved the manuscript.


\begin{thebibliography}{10}

\bibitem{lund2007renewable}
Henrik Lund.
\newblock Renewable energy strategies for sustainable development.
\newblock {\em Energy}, 32(6):912--919, 2007.

\bibitem{georgilakis2008technical}
Pavlos~S Georgilakis.
\newblock Technical challenges associated with the integration of wind power
  into power systems.
\newblock {\em Renewable and Sustainable Energy Reviews}, 12(3):852--863, 2008.

\bibitem{boyle1997renewable}
Godfrey Boyle et~al.
\newblock {\em Renewable energy: power for a sustainable future}.
\newblock Taylor \& Francis, 1997.

\bibitem{monteiro2009wind}
C~Monteiro, R~Bessa, V~Miranda, A~Botterud, J~Wang, G~Conzelmann, et~al.
\newblock Wind power forecasting: state-of-the-art 2009.
\newblock Technical report, Argonne National Laboratory (ANL), 2009.

\bibitem{braff2016value}
William~A Braff, Joshua~M Mueller, and Jessika~E Trancik.
\newblock Value of storage technologies for wind and solar energy.
\newblock {\em Nature Climate Change}, 6(10):964--969, 2016.

\bibitem{weitemeyer2015integration}
Stefan Weitemeyer, David Kleinhans, Thomas Vogt, and Carsten Agert.
\newblock Integration of renewable energy sources in future power systems: The
  role of storage.
\newblock {\em Renewable Energy}, 75:14--20, 2015.

\bibitem{weitemeyer2017optimal}
Stefan Weitemeyer, David Kleinhans, Lars Siemer, and Carsten Agert.
\newblock Optimal combination of energy storages for prospective power supply
  systems based on renewable energy sources, (submitted).
\newblock {\em Journal of Energy Storage}, 2017.

\bibitem{jones2014renewable}
Lawrence~E Jones.
\newblock {\em Renewable energy integration: practical management of
  variability, uncertainty, and flexibility in power grids}.
\newblock Academic Press, 2014.

\bibitem{schneider2007stochastic}
Johannes Schneider and Scott Kirkpatrick.
\newblock {\em Stochastic optimization}.
\newblock Springer Science \& Business Media, 2007.

\bibitem{ben2009robust}
Aharon Ben-Tal, Laurent El~Ghaoui, and Arkadi Nemirovski.
\newblock {\em Robust optimization}.
\newblock Princeton University Press, 2009.

\bibitem{grune2011nonlinear}
Lars Gr{\"u}ne and J{\"u}rgen Pannek.
\newblock Nonlinear model predictive control.
\newblock In {\em Nonlinear Model Predictive Control}, pages 43--66. Springer,
  2011.

\bibitem{siemer2016jes}
Lars Siemer, Frank Sch{\"o}pfer, and David Kleinhans.
\newblock Cost-optimal operation of energy storage units: Benefits of a
  problem-specific approach.
\newblock {\em Journal of Energy Storage}, 2016.

\bibitem{kleinhans2012estimation}
David Kleinhans.
\newblock Estimation of drift and diffusion functions from time series data: A
  maximum likelihood framework.
\newblock {\em Physical Review E}, 85(2):026705, 2012.

\bibitem{doob1942brownian}
Joseph~L Doob.
\newblock The brownian movement and stochastic equations.
\newblock {\em Annals of Mathematics}, pages 351--369, 1942.

\bibitem{gillespie1996exact}
Daniel~T Gillespie.
\newblock Exact numerical simulation of the ornstein-uhlenbeck process and its
  integral.
\newblock {\em Physical review E}, 54(2):2084, 1996.

\bibitem{maronna2006robust}
R.~Maronna, Douglas M., and Victor Y.
\newblock {\em Robust statistics}.
\newblock Wiley, 2006.

\bibitem{aissi2009min}
Hassene Aissi, Cristina Bazgan, and Daniel Vanderpooten.
\newblock Min--max and min--max regret versions of combinatorial optimization
  problems: A survey.
\newblock {\em European journal of operational research}, 197(2):427--438,
  2009.

\end{thebibliography}
\end{document}